\documentclass[graybox]{svmult}

\usepackage{mathptmx}       
\usepackage{helvet}         
\usepackage{courier}        
\usepackage{type1cm}        
\usepackage{makeidx}         
\usepackage{graphicx}        
\usepackage{multicol}        
\usepackage[bottom]{footmisc}

\makeindex             

\usepackage{amsmath, amssymb}
\usepackage{amsfonts,url}
\usepackage{hyperref}%
\usepackage{epsfig,tikz}

\newcommand\beq{\begin{equation}}
\newcommand\eeq{\end{equation}}
\newcommand\bce{\begin{center}}
\newcommand\ece{\end{center}}
\newcommand\bea{\begin{eqnarray}}
\newcommand\eea{\end{eqnarray}}
\newcommand\ba{\begin{array}}
\newcommand\ea{\end{array}}
\newcommand\ben{\begin{enumerate}}
\newcommand\een{\end{enumerate}}
\newcommand\bit{\begin{itemize}}
\newcommand\eit{\end{itemize}}
\newcommand\brr{\begin{array}}
\newcommand\err{\end{array}} 
\newcommand\bt{\begin{tabular}}  
\newcommand\et{\end{tabular}}

\newcommand\ul{\underline}
\newcommand\ol{\overline}

\renewcommand\S{{\mathcal S}}
\def\red{\operatorname{st}}
\def\inv{\operatorname{inv}}
\def\des{\operatorname{des}}
\def\O{O}

\def\N{{\mathcal N}}
\newcommand\Al{\operatorname{Allow}}
\newcommand\B{\operatorname{MF}}
\providecommand{\fr}[1]{\{#1\}}

\begin{document}

\title*{A survey of consecutive patterns in permutations}
\author{Sergi Elizalde}
\institute{Sergi Elizalde \at Department of Mathematics, Dartmouth College, Hanover, NH 03755, \email{sergi.elizalde@dartmouth.edu}}
\maketitle


\abstract{A consecutive pattern in a permutation $\pi$ is another permutation $\sigma$ determined by the relative order of a subsequence of contiguous entries of $\pi$. Traditional notions such as descents, runs and peaks can be viewed as particular examples of consecutive patterns in permutations, but the systematic study of these patterns has flourished in the last 15 years, during which a variety of different techniques have been used.
We survey some interesting developments in the subject, focusing on exact and asymptotic enumeration results, the classification of consecutive patterns into equivalence classes, and their applications to the study of one-dimensional dynamical systems.
}

\section{Introduction}

Patterns in permutations are implicit already in the work of MacMahon from one century ago, but 
interest in them has grown in the last few decades, inspired by work of Knuth~\cite{Knu} on sorting permutations using data structures such as stacks and double-ended queues. He showed that permutations sortable with these devices
can be characterized as those not containing subsequences of length 3 whose entries are in a prescribed order (we call such a subsequence a classical pattern), and that they are counted by the Catalan numbers.
A few years later, the systematic enumeration of permutations avoiding classical patterns was continued by Simion and Schmidt~\cite{SS}. Interest in permutation patterns has grow steadily since then, becoming a very active area of research, with connections with other areas of mathematics, such as algebraic geometry, and other fields such as computer science, biology, and physics.

In this survey we are interested in an important variation to the definition of permutation patterns used in~\cite{Knu,SS}, namely the notion of {\em consecutive} patterns. 
In an occurrence of a consecutive pattern in a permutation, the positions of the entries are required to be adjacent.
Long before this concept was introduced and systematically studied in~\cite{EliNoy} by Elizalde and Noy, there had been implicit appearances of consecutive patterns in the combinatorics literature. For instance,
occurrences of $21$ are descents, whose distribution is given by the Eulerian numbers;
occurrences of $321$ are sometimes called double descents~\cite{FV};
occurrences of $132$ and $231$ are called peaks and play an important role in algebraic combinatorics;
maximal occurrences of $12\dots m$ are called runs;
and permutations avoiding $123$ and $321$ are the well-known alternating permutations~\cite{And}, counted by the Euler numbers. 

In the last 15 years, a growing number of authors have made significant progress on the study of consecutive patterns in permutations by using an array of different techniques. Despite this ongoing research, many important questions in the field remain unanswered. In this article we survey some of these results and discuss how they are related, without the goal of being exhaustive.

After introducing the basic definitions and notation in Section~\ref{sec:definitions}, Section~\ref{sec:Wilf} discusses an equivalence relation on patterns which is analogous to Wilf equivalence for classical patterns. Section~\ref{sec:enumeration} surveys exact enumeration formulas that have been obtained using approaches such as the symbolic method, binary trees, the cluster method, symmetric funcions, homological algebra, and spectral theory. In Section~\ref{sec:asym} we summarize some results concerning the asymptotic behavior of the number of permutations avoiding a consecutive pattern. Finally, Section~\ref{sec:dynamical} presents a recently discovered application of consecutive patterns to the study of time series.

We will not discuss the extensive literature on classical (i.e. non-consecutive) patterns, for which we refer the reader to the surveys~\cite{Vatter} and~\cite{Kitaev}. We point out that even though some of the questions in both areas are similar, most of the techniques that have been used are significantly different, and so are the results that have been obtained, from the classification into equivalence classes to the
description of the asymptotic behavior of permutations avoiding a pattern. Exact enumeration generally seems easier in the consecutive case: the only infinite family of classical patterns for which a generating function is known is the monotone family, whereas for consecutive patterns the distribution of occurrences for several families of patterns is known. On the other hand, some small classical patterns give rise to simple algebraic generating functions, which does not happen in the consecutive case.

\section{Definitions}\label{sec:definitions}

Let $\S_n$ denote the set of permutations of $\{1,2,\dots,n\}$.
For a sequence $\tau$ of $k$ distinct positive integers, let $\red({\tau})$ denote the permutation in $\S_k$ obtained by
replacing the smallest entry of $\tau$ with~$1$, the second smallest with~$2$, and so on. For example, $\red(394176) = 263154$.
Given $\pi=\pi_1\pi_2\cdots\pi_n\in\S_n$ and $\sigma=\sigma_1\sigma_2\cdots\sigma_k\in\S_k$,
an occurrence of $\sigma$ in $\pi$ as a \emph{consecutive pattern} is a subsequence of $k$ contiguous entries of $\pi$ such that $\red(\pi_{i}\pi_{i+1}\cdots\pi_{i+k-1})=\sigma_1\sigma_2\cdots\sigma_k$.
For example, in $\pi=132546$, the subsequences $1325$ and $2546$ are two occurrences of the pattern $\sigma=1324$.

Denote by  $c_\sigma(\pi)$ the number of occurrences of $\sigma$ in $\pi$ as a consecutive pattern. If $c_\sigma(\pi)=0$, we say that $\pi$ {\em avoids} $\sigma$. For example, $25134$ avoids $132$. Let $\alpha_n(\sigma)$
be the number of permutations in $\S_n$ that avoid $\sigma$ as a consecutive pattern. The notions of occurrence, containment and avoidance in this survey will refer to consecutive patterns unless explicitly stated otherwise, in which case we will use the term {\em classical pattern} to mean that 
the entries in an occurrence have no adjacency requirement. 
 
Let
$$P_{\sigma}(u,z)=\sum_{n\ge0} \sum_{\pi\in\S_n} u^{c_\sigma(\pi)}\frac{z^n}{n!}$$ be the exponential generating function for occurrences of $\sigma$.
Setting $u=0$, we get the generating function for $\sigma$-avoiding permutations, $P_\sigma(0,z)=\sum_{n\ge0} \alpha_n(\sigma)\frac{z^n}{n!}$.

A useful invariant of a pattern $\sigma$ is its overlap set $\O_\sigma$, defined  as the set of indices $i$ with $1\le i<m$
such that $\red(\sigma_{i+1}\sigma_{i+2}\dots\sigma_m)=\red(\sigma_1\sigma_2\dots \sigma_{m-i})$.
Equivalently, $i\in\O_\sigma$ if two occurrences of $\sigma$ in a permutation can have starting positions at distance $i$ from each other.
For example $O_{1423}=\{2,3\}$.
Note that $m-1\in \O_\sigma$ for every $\sigma\in\S_m$. If $m\ge3$, a pattern $\sigma\in\S_m$ for which $\O_\sigma=\{m-1\}$ is said to be
{\em non-overlapping} (sometimes also called {\em minimally overlapping} or {\em non-self-overlapping}). Equivalently, $\sigma$ is non-overlapping if two occurrences of $\sigma$ in a permutation cannot overlap in more than one position.
For example, the patterns $132$, $1243$, $1342$, $21534$ and $34671285$ are non-overlapping. We denote by $\N_m$ the set of non-overlapping patterns in $\S_m$.
Non-overlapping patterns have been studied by Duane and Remmel~\cite{DR} and by B\'ona~\cite{Bon},
who showed that $\lim_{m\to\infty}|\N_m|/m!\approx 0.364$.

\section{Consecutive Wilf-equivalence classes}\label{sec:Wilf}

In analogy with the definition of Wilf-equivalence for classical patterns, we say that two patterns $\sigma$ and $\tau$ are {\em c-Wilf-equivalent} if $P_{\sigma}(0,z)=P_{\tau}(0,z)$ (that is,
$\alpha_n(\sigma)=\alpha_n(\tau)$ for all $n$), and we say that they are {\em strongly c-Wilf-equivalent} if $P_\sigma(u,z)=P_\tau(u,z)$. 
For example, every pattern $\sigma=\sigma_1 \cdots \sigma_m$ is strongly c-Wilf-equivalent to its reversal $\sigma_m \cdots \sigma_1$ and its complement $(m{+}1{-}\sigma_1) \cdots (m{+}1{-}\sigma_m)$. Some  non-trivial examples of strong c-Wilf-equivalent pairs are 
$1342\sim1432$ and $154263\sim165243$.
An open problem in consecutive patterns (whose analogue in the classical case is wide open as well) is to classify patterns into these equivalence classes. 
It was conjectured by Nakamura that the word {\em strongly} in the above definition is superfluous, in sharp contrast with the case of classical patterns:

\begin{conjecture}[{\cite[Conjecture 5.6]{Nak}}]\label{conj:Nak}
If two patterns are c-Wilf-equivalent then they are also strongly c-Wilf-equivalent.
\end{conjecture}

This conjecture is known to hold in the special case of non-overlapping patterns (see~\cite[Lemma 3.2]{EliCMP} and \cite[Theorem 11]{MR}). In fact, 
Dwyer and Elizalde~\cite{DE} have recently proved that if two non-overlapping patterns $\sigma$ and $\tau$ are c-Wilf-equivalent, then 
they are {\em super-strongly} c-Wilf equivalent, meaning that for every $n$ and for every set of positions $S$, the number of permutations in $\S_n$ having occurrences of $\sigma$ exactly at positions $S$ is the same as for $\tau$. For example, $1342$ and $2431$ are super-strongly c-Wilf equivalent, but $1423$ and $3241$ are not.
It is worth mentioning that Conjecture~\ref{conj:Nak} has an analogue for words \cite[Conjectiure 1.2]{PV}, which, as Pantone and Vatter showed, would imply the so-called Rearrangement Conjecture from~\cite{KLRS}.

So far, (strong) c-Wilf-equivalence classes have been characterized for patterns of length up to 6. 
There are $2$ equivalence classes of patterns of length~$3$, represented by the patterns $123$ and $132$; $7$ classes of patterns of length~$4$, represented by $1234$, $2413$, $2143$, $1324$, $1423$, $1342$, and $1243$ (see \cite{EliNoy}); $25$ classes for patterns of length~$5$; and $92$ for patterns of length~$6$ (see \cite{Nak,EliNoy2}).

For non-overlapping patterns of any given length $m$, it was shown in~\cite{DK,DR} ---proving a conjecture from~\cite{Elitesis}--- that the first and last entry of a pattern determine the strong c-Wilf-equivalence class where it belongs. The converse, once we account for the fact that the reversal and the complement of a pattern preserve its c-Wilf-equivalence class, is conjectured to be true in~\cite{EliCMP}.
The above result has been generalized to sets of (possibly overlapping) consecutive patterns by Khoroshkin and Shapiro~\cite{KS}, by giving a list of sufficient
conditions for such sets to be strongly c-Wilf-equivalent. As shown in~\cite{DE}, the same conditions imply super-strong c-Wilf-equivalence.

\section{Exact enumeration}\label{sec:enumeration}

\subsection{The origins}

As mentioned in the introduction, some traditional results on permutations can be interpreted as statements about consecutive patterns, such as the well-known generating function for the Eulerian numbers,
$$P_{21}(u,z)=\frac{1-u}{e^{(u-1)z}-u},$$
and Andr\'e's generating function $\sec z+\tan z$ for the Euler numbers, which count up-down permutations~\cite{And}.

Almost one century later, the first published formula for the number of
permutations with no increasing runs of length $m$ or more appears in a 1962 book of David and Barton~\cite[pg. 156--157]{DB}:
\begin{equation}\label{eq:DB}
P_{12\dots m}(0,z)=\left(\sum_{j\ge0}\frac{z^{jm}}{(jm)!}-\sum_{j\ge0}\frac{z^{jm+1}}{(jm+1)!}\right)^{-1}.
\end{equation}
Soon after, we find Entringer's enumeration of permutations according to their number of cyclic peaks (which he called maxima)~\cite{Ent}, and Carlitz and Scoville's enumeration of permutations according to rising maxima and falling maxima~\cite{CS}, which are essentially (up to the definition at the beginning and at the end of the permutation) occurrences of $132$ and $231$.

It was not until 2001 that the systematic study of consecutive patterns in permutations was tackled. In~\cite{EliNoy}, Elizalde and Noy introduced the notion of consecutive patterns and found expressions for the generating functions $P_\sigma(u,z)$ for certain families of patterns of arbitrary length, which include those of length 3 and three of the seven c-Wilf equivalence classes of patterns of length 4.
Letting $$\omega_\sigma(u,z)=\frac{1}{P_\sigma(u,z)},$$ they gave differential equations satisfied by $\omega_\sigma(u,z)$.
The proofs are based on representations of permutations as increasing binary trees, together with the symbolic method of Flajolet and Sedgewick~
\cite{FS}. In the following statements, the derivatives are taken with respect to $z$.

\begin{theorem}[\cite{EliNoy}]\label{thm:monotone}
Let $\sigma=12\dots m$, where $m\ge3$. Then $\omega=\omega_\sigma(u,z)$ is the solution of
\beq\label{eq:omegamon}\omega^{(m-1)}+(1-u)(\omega^{(m-2)}+\dots+\omega'+\omega)=0\eeq
with $\omega(0)=1,\omega'(0)=-1,\omega^{(i)}(0)=0$ for $2\le i\le m-2$.
\end{theorem}

When $u=0$, the solution of the above differential equation can be expressed as a series, recovering Equation~\eqref{eq:DB}.

Another family enumerated in~\cite{EliNoy} consists of certain non-overlapping patterns that begin with $1$. The above observation that the first and last entry determine the (strongly) c-Wilf-equivalence class of a non-overlapping pattern allows us to state the theorem more generally as follows.

\begin{theorem}[\cite{EliNoy}]\label{thm:1b}
Let $\sigma\in\N_m$ with $\sigma_1=1$, and let $b=\sigma_m$. Then $\omega=\omega_\sigma(u,z)$ is the solution of
\beq\label{eq:omega1b}\omega^{(b)}+(1-u)\frac{z^{m-b}}{(m-b)!}\omega'=0\eeq
with $\omega(0)=1,\omega'(0)=-1,\omega^{(i)}(0)=0$ for $2\le i\le b-1$.
\end{theorem}

When $b=2$, the above differential equation can be solved to obtain
$$P_\sigma(u,z)=\left(1-\int_0^z e^{(u-1)t^{m-1}/(m-1)!} dt\right)^{-1}$$
for any $\sigma\in\S_m$ with $\sigma_1=1$ and $\sigma_m=2$ (note that such a pattern is always non-overlapping).

A related problem considered in~\cite{EliNoy} is the enumeration of permutations that simultaneously avoid several consecutive patterns. Many more results in this direction, which we do not include here, have been obtained by Kitaev~\cite{Kit03}, Kitaev and Mansour~\cite{KM,KM2}, and Aldred, Atkinson and McCaughan~\cite{AAM}.

In~\cite{Kit05}, Kitaev considered the maximum number of non-ovelapping occurrences (i.e., not sharing any entries) of a pattern $\sigma$ in a permutation $\pi$, which we denote by $\mathrm{nlap}_\sigma(\pi)$. He showed that the distribution of this statistic depends only on the number of $\sigma$-avoiding permutations:

\begin{theorem}[\cite{Kit05}]\label{thm:kitaev} For every $\sigma$,
$$\sum_{n\ge0} \sum_{\pi\in\S_n} u^{\mathrm{nlap}_\sigma(\pi)}\frac{z^n}{n!}=\frac{P_\sigma(0,z)}{(1-u)+u(1-z)P_\sigma(0,z)}.$$
\end{theorem}

\subsection{The cluster method}
A powerful method to obtain expressions for $P_\sigma(u,z)$ is the {\em cluster method} of Goulden and Jackson~\cite{GJ79}, which is based on inclusion-exclusion (see also \cite[Ch. 2.8]{GJ}, its restatement in \cite[Ex. 4.40]{EC1}, and its extensions in~\cite{NZ}). Originally used for counting words according to occurrences of factors (i.e. consecutive subwords), it can be formulated in terms of permutations as stating that
\beq\label{eq:clustermethod} P_\sigma(u,z)=\frac{1}{1-z-R_\sigma(u-1,z)},\eeq
where $R_\sigma(t,z)$ is an exponential generating function where the coefficient of $t^kz^n/n!$ is the number of $k$-clusters with respect to $\sigma$, that is, permutations of length $n$ with $k$ marked occurrences of $\sigma$ with the property that each marked occurrence overlaps with the next one, and the whole permutation is covered with these marked occurrences. Note that a cluster can have other occurrences of $\sigma$ that are not marked. For example, $\ul{1\ 6\ \ol{2\ 8\ 3}}\ol{\ 11 \ }\ul{\ol{4}\ 9\ 5\ 10\ 7}$ is a $3$-cluster with respect to $14253$, where the horizontal lines indicate the three marked occurrences.

To prove Equation~\eqref{eq:clustermethod}, note that the generating function $1/(1-z-R_\sigma(t,z))$ counts permutations obtained as sequences of single entries and clusters, where the exponent of $t$ is the number of marked occurrences of $\sigma$ in these clusters. Thus, a permutation of length $n$ with $\ell$ occurrences of $\sigma$ contributes a term $t^k z^n/n!$ for each of the $\binom{\ell}{k}$ ways to mark $k$ of these $\ell$ occurrences, for every $k\le\ell$.
Making the substitution $t=u-1$, the contribution of such a permutation to the right hand side of~\eqref{eq:clustermethod} is
$$\sum_{k=0}^\ell \binom{\ell}{k} (u-1)^k \frac{z^n}{n!}=u^\ell \frac{z^n}{n!},$$
and thus equal to the contribution to the left hand side.

Equation~\eqref{eq:clustermethod} is useful because it reduces the enumeration of occurrences of $\sigma$ in all permutations to the enumeration of clusters, which are typically more tractable. For example, when $\sigma=12\dots m$, clusters are simply increasing permutations with marked overlapping occurrences of $\sigma$. When $u=0$, Equation~\eqref{eq:DB} follows from this observation with a little work. When $\sigma$ is a non-overlapping pattern as in the statement of Theorem~\ref{thm:1b}, computing the number of clusters and applying Equation~\eqref{eq:clustermethod} we get
\beq\label{eq:cluster1b} P_\sigma(u,z)=\left(1-z-\sum_{k\ge1}\prod_{j=2}^k\binom{j(m-1)-b+1}{m-b}\frac{(u-1)^k z^{(m-1)k+1}}{((m-1)k+1)!}\right)^{-1}.\eeq
In~\cite{Kit05}, Kitaev used inclusion-exclusion arguments to deduce Equation~\eqref{eq:cluster1b} when $u=0$.

Rawlings~\cite{Raw} noticed the power of the cluster method and used it to obtain enumeration formulas for the patterns $12\dots m$, 
$12\dots(m-2)m(m-1)$ and $1m(m-1)\dots2$ with an additional variable that keeps track of the number of inversions in the permutation, thus generalizing some of the results in~\cite{EliNoy}. Although his proof relies on a bijection of F\'edou, it is in fact possible to directly generalize the above argument proving Equation~\eqref{eq:clustermethod} by adding a variable $q$ that marks the number of inversions (both in the permutations on the left hand side and on the clusters on the right hand side), and using generating functions where the $n!$ in the exponential generating functions is replaced by its $q$-analogue $[n]_q!=\prod_{i=1}^n(1+q+\dots+q^{i-1})$. 

For example, refining Equation~\eqref{eq:cluster1b} by applying this $q$-analogue of the cluster method, we get the following generalization of 
Rawlings' formulas, where $\inv(\pi)$ denotes the number of inversions of $\pi$, and $\binom{n}{k}_q=\frac{[n]_q!}{[k]_q![n-k]_q!}$ is the $q$-binomial coefficient:
\begin{theorem}
Let $\sigma\in\N_m$ with $\sigma_1=1$, and let $b=\sigma_m$. Then 
\begin{multline*}
\sum_{n\ge0} \sum_{\pi\in\S_n} q^{\inv(\pi)}u^{c_\sigma(\pi)}\frac{z^n}{[n]_q!}\\
=\left(1-z-\sum_{k\ge1}\prod_{j=2}^k\binom{j(m-1)-b+1}{m-b}_q\frac{q^{k\inv(\sigma)}(u-1)^k z^{(m-1)k+1}}{[(m-1)k+1]_q!}\right)^{-1}.\end{multline*}
\end{theorem}

Similarly, denoting by $\S_n(12\dots m)$ the set of permutations in $\S_n$ that avoid $12\dots m$, 
the $q$-analogue of Equation~\eqref{eq:DB} is
$$\sum_{n\ge0} \sum_{\pi\in\S_n(12\dots m)} q^{\inv(\pi)}\frac{z^n}{[n]_q!}=\left(\sum_{j\ge0}\frac{z^{jm}}{[jm]_q!}-\sum_{j\ge0}\frac{z^{jm+1}}{[jm+1]_q!}\right)^{-1}.$$

Dotsenko and Khoroshkin~\cite{DK} developed a variation of the cluster method that they call the chain method. What makes their approach particularly interesting is that it is motivated by homological algebra, and more specifically by the computation of free resolutions of Anick type for shuffle algebras. As it turns out, shuffle algebras with monomial relations have bases that can be naturally described via generalized colored permutations avoiding consecutive patterns. From a combinatorial perspective, the chain method is also based on inclusion-exclusion, and the role of clusters is played by chains, which are clusters with additional restrictions. This reduces the number of objects to consider by eliminating some of the clusters whose contributions to $R_\sigma(-1,z)$ cancel out (for example, permutations that are simultaneously a $k$-cluster and a $k+1$-cluster). Recurrences for the number of clusters and chains are given in~\cite{DK} for some patterns.
Further applications of this homological algebra approach to consecutive patterns are studied by Khoroshkin and Shapiro in~\cite{KS},
and generalizations to patterns in labeled trees are considered by Dotsenko~\cite{Dottrees}.

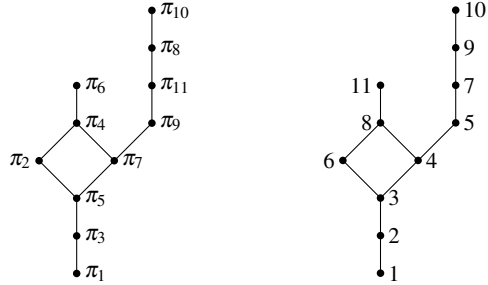
\begin{figure}[htb]
\centering
\begin{tikzpicture}[scale=.5]
\fill (1,0) circle (0.1) node[right]{$\pi_1$};
\fill (1,1) circle (0.1) node[right]{$\pi_3$};
\fill (1,2) circle (0.1) node[right]{$\pi_5$};
\fill (0,3) circle (0.1) node[left]{$\pi_2$};
\fill (2,3) circle (0.1) node[right]{$\pi_7$};
\fill (1,4) circle (0.1) node[right]{$\pi_4$};
\fill (1,5) circle (0.1) node[right]{$\pi_6$};
\fill (3,4) circle (0.1) node[right]{$\pi_9$};
\fill (3,5) circle (0.1) node[right]{$\pi_{11}$};
\fill (3,6) circle (0.1) node[right]{$\pi_8$};
\fill (3,7) circle (0.1) node[right]{$\pi_{10}$};
\draw (1,0)--(1,2)--(0,3)--(1,4)--(1,5);
\draw (1,2)--(3,4)--(3,7);
\draw (2,3)--(1,4);
\end{tikzpicture}\hspace{15mm}
\begin{tikzpicture}[scale=.5]
\fill (1,0) circle (0.1) node[right]{$1$};
\fill (1,1) circle (0.1) node[right]{$2$};
\fill (1,2) circle (0.1) node[right]{$3$};
\fill (0,3) circle (0.1) node[left]{$6$};
\fill (2,3) circle (0.1) node[right]{$4$};
\fill (1,4) circle (0.1) node[left]{$8$};
\fill (1,5) circle (0.1) node[left]{$11$};
\fill (3,4) circle (0.1) node[right]{$5$};
\fill (3,5) circle (0.1) node[right]{$7$};
\fill (3,6) circle (0.1) node[right]{$9$};
\fill (3,7) circle (0.1) node[right]{$10$};
\draw (1,0)--(1,2)--(0,3)--(1,4)--(1,5);
\draw (1,2)--(3,4)--(3,7);
\draw (2,3)--(1,4);
\end{tikzpicture}
\caption{A cluster poset for $\sigma=14253$ and a linear extension corresponding to the $3$-cluster $\ul{1\ 6\ \ol{2\ 8\ 3}}\ol{\ 11 \ }\ul{\ol{4}\ 9\ 5\ 10\ 7}$.}
\label{fig:lext}
\end{figure}

In~\cite{EliNoy2}, Elizalde and Noy observed that counting clusters is equivalent to counting linear extensions of certain posets, which they call {\em cluster posets}. For example, the cluster posets for the monotone pattern are chains, and the cluster posets for non-overlapping patterns have a simple tree-like structure. Another example is the poset describing the order relations satisfied by the entries of a cluster $\pi\in\S_{11}$ with marked occurrences of $\sigma=14253$ starting at positions $1,3,7$, which is shown in Figure~\ref{fig:lext}.
The interpretation in terms of linear extensions is exploited in~\cite{EliNoy2} to obtain differential equations satisfied by $\omega_\sigma(u,z)$ for several patterns. For example, for all patterns for which the cluster posets are chains (this includes $12435$ and $124536$, for instance), an equation similar to Theorem~\ref{thm:monotone} is derived, where the derivatives that appear depend on the overlap set $\O_\sigma$. Differential equations are given for other patterns of length $5$ and $6$ as well, and for the pattern $1324$, which has also been studied in~\cite{DK} and~\cite{LR}.

\subsection{Symmetric functions and brick tabloids}

Before the cluster method became a common tool in the study of consecutive patterns, Mendes and Remmel~\cite{MR} developed in 2006 a remarkable technique to express generating functions for permutations and words avoiding consecutive patterns in terms of the reciprocal of another power series.
Their method combines symmetric function identities expressed in terms of what they call {\em brick tabloids}, ring homomorphisms from symmetric functions to polynomials (inspired by Brenti's work~\cite{Brenti}), and sign-reversing involutions on brick tabloids. 

The first application in~\cite{MR} is an alternative proof of Theorem~\ref{thm:kitaev}, which is generalized to avoidance of sets of patterns, and refined by the inversion number.
Another important result from \cite{MR} is an expression of $P_\sigma(0,z)$ as the reciprocal of a power series whose coefficients are signed sums of permutations having occurrences of $\sigma$ at prescribed positions determined by a certain set associated to $\sigma$. This power series is reminiscent of the cluster generating function $R_\sigma(-1,z)$, and in the case of non-overlapping patterns it actually coincides with it~\cite{DR}. Subsequently, Liese and Remmel~\cite{LR} gave explicit expressions for this power series in some cases where $\sigma$ is a shuffle of an increasing sequence with another pattern.

One advantage of the Mendes--Remmel method is that it often yields refinements by statistics such as the number of descents, the number of inversions or the major index, and it extends naturally to permutations avoiding sets of patterns, to tuples of permutations, and to colored permutations (elements of the wreath product $C_k\wr\S_n$).
An example of such a refinement is the following, where $\des(\pi)$ denotes the number of descents of $\pi$.
\begin{theorem}[\cite{MR}] For $m\ge2$,
$$\sum_{n\ge0} \sum_{\pi\in\S_n(m\dots21)} x^{\des(\pi)} \frac{z^n}{n!}=\left(\sum_{n\ge0}\sum_{i\ge0}(-1)^{n-i} r_{n,i,m}\,x^i\frac{z^n}{n!}\right)^{-1},$$
where $r_{n,i,m}$ is the number of rearrangements of $i$ zeroes and $n-1-i$ ones such that $m-1$ zeroes never appear consecutively.
\end{theorem}
Setting $x=1$, Equation~\eqref{eq:DB} follows by applying a sign-reversing involution to the above rearrangements.

On the other hand, most of the formulas resulting from this method do not keep track of the number of occurrences ({\em matches}, in the terminology from~\cite{MR}) of the pattern, unlike with the cluster method. Nonetheless, there is a special type of patterns for which Mendes and Remmel~\cite[Theorem 11]{MR} express $P_\sigma(u,z)$ in terms of $P_\sigma(0,z)$. These are precisely the non-overlapping patterns, as shown in~\cite{DE}. We remark that for $\sigma\in\N_m$, the cluster method can also be used to show that
$$\omega_\sigma(u,z)=1+\frac{\omega_\sigma(0,(1-u)^{\frac{1}{m-1}}\,z)-1}{(1-u)^{\frac{1}{m-1}}}.$$

\subsection{Other questions and approaches}\label{sec:other}

Even though formulas for $P_\sigma(u,z)$ are known for several patterns, the nature of these generating functions for an arbitrary $\sigma$ is not well understood. For the case of patterns without the adjacency requirement, it was conjectured by Noonan and Zeilberger~\cite{NZ96}, extending a speculation of Gessel~\cite{Ges}, that the generating function for permutations avoiding a classical pattern is always D-finite, that is, it satisfies a linear differential equation with polynomial coefficients. This conjecture has very recently been disproved in~\cite{GP}.
In the case of consecutive patterns, the interesting question is not whether $P_\sigma(0,z)$ is D-finite (since it clearly is not, already for $\sigma=123$) but whether its reciprocal $\omega_\sigma(0,z)$ is D-finite. 
Even though it is for the patterns for which the enumeration is known, Elizalde and Noy conjecture that this is not the case in general, and they suggest that $1423$ is a counterexample.

\begin{conjecture}[\cite{EliNoy2}]
The generating function $\omega_{1423}(0,z)$ is not D-finite.
\end{conjecture}

An automated approach to the enumeration of permutations avoiding consecutive patterns has been carried out by Baxter, Nakamura and Zeilberger~\cite{BNZ}. Their algorithms, implemented in Maple, automatically derive functional equations for the generating functions, and use them to generate the corresponding sequences in polynomial time. In related work, Baxter and Pudwell~\cite{BP} provide an algorithm that, given any set of consecutive patterns, returns a recurrence (called an {\em enumeration scheme}) that can be used to efficiently compute the number of permutations avoiding them.

In a different direction, some recent papers study the infinite poset of all permutations $\cup_{n\ge1}\S_n$ where one defines
$\sigma\le\pi$ if $\pi$ contains $\sigma$ as a consecutive pattern. In~\cite{BFS}, Bernini, Ferrari and Steingr\'{\i}msson gave a recursive formula for the M\"obius function of this poset. In~\cite{SW}, Sagan and Willenbring used discrete Morse theory to determine its homotopy type. Further properties of this poset are currently being studied in~\cite{EliMc}.

In recent years, a number of variations and generalizations of consecutive and classical patterns have been studied in the literature. These include {\em vincular} patterns, originally called generalized patterns when first defined by Babson and Steingr\'{\i}msson~\cite{BS}; {\em bivincular} patterns, introduced by Bousquet-M\'elou et al.~\cite{BCDK}; {\em mesh} patterns, introduced by Br\"and\'en and Claesson~\cite{BC}; and {\em partially ordered} patterns, proposed by Kitaev~\cite{Kit05,Kit10}.

\section{Asymptotic enumeration}\label{sec:asym}

\subsection{Growth rates}
Even in cases where no formula or generating function for the numbers $\alpha_n(\sigma)$ is known, it is natural to ask how fast this sequence grows. The first general result in this direction is the following.

\begin{proposition}[\cite{Eliasym}]\label{prop:lim} For every $\sigma\in\S_m$ with $m\ge3$,
the limit $$\lim_{n\rightarrow\infty}\left(\frac{\alpha_n(\sigma)}{n!}\right)^{1/n}$$ exists, and it is strictly between $0$ and $1$.
\end{proposition}

This limit, which we denote by $\rho_\sigma$, is called the {\em growth rate} of $\sigma$. It is a well-known fact in analytic combinatorics~\cite[Theorem IV.7]{FS} that $\rho_\sigma^{-1}$ equals the the radius of convergence of $P_\sigma(0,z)$ as a function of complex variable, and that $P_\sigma(0,z)$ has a real singularity at $z=\rho_\sigma^{-1}$ \cite[Theorem IV.6]{FS} and no singularities in $|z|<\rho_\sigma^{-1}$.

The constants $\rho_\sigma$ were computed in~\cite{EliNoy} for some small patterns.
For example, $\rho_{123}=\frac{3\sqrt{3}}{2\pi}\approx 0.8269933$ and $\rho_{132}\approx 0.7839769$, the reciprocal of the unique positive root of $\int_0^z e^{-t^2/2} dt=1$. For monotone patterns, it can be shown using Equation~\eqref{eq:DB} and the above discussion that, letting $m$ go to infinity,
$$\rho_{12\dots m}=1-\frac{1}{m!}+\frac{1}{(m+1)!}+O\left(\frac{1}{m!^2}\right).$$

In general, approximate values for $\rho_\sigma$ can now be found quickly using the Maple packages from~\cite{BNZ}. 
For the patterns of length 3 and 4 solved in~\cite{EliNoy}, it was also shown that as $n$ tends to infinity, $\alpha_n(\sigma)\sim\gamma_\sigma\,\rho_\sigma^n\,n!$ for some constant $\gamma_\sigma$. Warlimont~\cite{War} conjectured that this asymptotic behavior holds for all patterns $\sigma$, and Ehrenborg, Kitaev and Perry~\cite{EKP} proved this fact using methods from spectral theory of integral operators:

\begin{theorem}[\cite{EKP}]\label{thm:EKP}
For every $\sigma$, $$\frac{\alpha_n(\sigma)}{n!}=\gamma_\sigma\rho_\sigma^n+O(r_\sigma^n)$$ as $n$ tends to infinity, for some constants $\gamma_\sigma$ and $r_\sigma<\rho_\sigma$.
\end{theorem}

The techniques in~\cite{EKP} give detailed asymptotic expansions. An alternative proof of Theorem~\ref{thm:EKP} seems to be within reach of
the standard methods from singularity analysis of generating functions. It was shown in \cite{EliCMP} that $\rho_\sigma^{-1}$ is the smallest positive zero of $\omega_\sigma(0,z)$, and that this zero is simple. Proving that this is the unique zero of $\omega_\sigma(0,z)$ of minimum modulus would imply Theorem~\ref{thm:EKP} by~\cite[Theorem IV.10]{FS}.

It is not known whether $\omega_\sigma(0,z)$ is always an entire function in the complex plane. It was shown in~\cite{EliNoy2} that it is entire for a large family of patterns, but that $\omega_{2413}(2,z)$ is not.

\subsection{Comparisons between patterns}

Another interesting direction of research aims at comparing the growth rates of different patterns. Phrased in a different way, one question is to determine, given two patterns $\sigma$ and $\tau$, whether $\alpha_n(\sigma)>\alpha_n(\tau)$ for large enough $n$. Two instances of this question were answered already in~\cite{EliNoy}, where it is shown that 
\begin{equation}\label{eq:compare}
\alpha_n(123)>\alpha_n(132) \mbox{ for }n\ge4, \quad \mbox{and}\quad \alpha_n(1342)>\alpha_n(1243) \mbox{ for }n\ge7.
\end{equation}

Even though the question remains wide open in general, the following result describes which patterns of any given length maximize and minimize $\rho_\sigma$.

\begin{theorem}[\cite{EliCMP}]\label{thm:CMP}
For every $\sigma\in\S_m$ there exists $n_0$ such that $$\alpha_n(12\dots(m{-}2)m(m{-}1)) \le \alpha_n(\sigma)\le \alpha_n(12\dots m)$$ for all $n\ge n_0$.
\end{theorem}

It is shown in \cite{EliCMP} that the above inequalities are strict unless $\sigma$ is a trivial symmetry of one of the above two patterns. The fact that $\alpha_n(\sigma)\le \alpha_n(12\dots m)$ for all $\sigma\in\S_m$ had been conjectured in~\cite{EliNoy}, and proved for $\sigma\in\N_m$ in~\cite{EliNoy2}. The fact that $\alpha_n(12\dots(m{-}2)m(m{-}1)) \le \alpha_n(\sigma)$ for all $\sigma\in\S_m$  had been conjectured in~\cite{Nak}.

It is interesting to note that, perhaps counter-intuitively, the analogous statements for classical patterns do not hold, as shown by B\'ona~\cite{Bon97}. In fact, $12\dots m$ and $12\dots (m{-}2)m(m{-}1)$ are Wilf-equivalent as classical patterns.
A different analogue of the upper bound for another definition of pattern avoidance, which requires adjacent positions and adjacent values, was considered in~\cite{Bon08}.

The proof of Theorem~\ref{thm:CMP} combines singularity analysis of generating functions with combinatorial arguments involving linear extensions of cluster posets. The inequality $\alpha_n(\sigma)<\alpha_n(12\dots m)$ is proved by showing that $\rho_\sigma<\rho_{12\dots m}$ if $\sigma$ is not monotone. After arguing that $\rho_\sigma^{-1}$ is the smallest zero of $\omega_\sigma(0,z)$, the proof can be reduced to showing that
$\omega_{12\dots m}(0,z)<\omega_\sigma(0,z)$ for $0<z<C$, for a certain constant $C$.
This is then done by carefully bounding the terms of the cluster generating functions $R_\sigma(-1,z)$.

A version of Theorem~\ref{thm:CMP} for non-overlapping patterns also appears in~\cite{EliCMP}. It 
states that for every $\sigma\in\N_m$, there exists $n_0$ such that 
\begin{equation}\label{eq:bound_nol}
\alpha_n(12\dots(m{-}2)m(m{-}1))\le\alpha_n(\sigma)\le\alpha_n(134\dots m2)
\end{equation}
for all $n\ge n_0$.

A probabilistic proof of the right inequality in Theorem~\ref{thm:CMP} for $m\ge 5$ has been found by Perarnau~\cite{Perarnau} using Suen's inequality, which gives an upper bound on the probability that none of the events in a certain collection occur simultaneously. The probabilistic method is suitable for consecutive patterns because the events that describe occurrences of $\sigma$ in a large permutation have few dependencies.
It is also shown in~\cite{Perarnau} that for $\sigma\in\S_m$, as $m$ goes to infinity, $$\rho_\sigma\ge 1-\frac{1}{m!}-O\left(\frac{m}{m!^2}\right),$$
and that the growth rate of most patterns is close to this lower bound.

We point out that even though the inequalities in Equation~\eqref{eq:compare} were proved combinatorially by explicitly constructing injective maps \cite{EliNoy}, no combinatorial proofs of Theorem~\ref{thm:CMP} and Equation~\eqref{eq:bound_nol} are known.

\section{Applications to dynamical systems}\label{sec:dynamical}

Aside from their interest from a purely combinatorial perspective, consecutive patterns arise naturally in the study of time series, where they play a role in distinguishing deterministic from random sequences (see~\cite{Amigobook} for a survey).
Given a map $f$ from a linearly ordered set (such as an interval in the real line) to itself, consider the finite sequences (orbits) that are obtained by iterating $f$, starting from different initial points: $x,f(x),f(f(x)),\dots,f^{n-1}(x)$. If these values are all different, their relative order determines a permutation $$\red(x,f(x),f(f(x)),\dots,f^{n-1}(x))\in\S_n.$$
Permutations obtained in this way, for varying $x$ in the domain of $f$, are called {\em allowed patterns} of $f$, or patterns {\em realized} by $f$, and denoted by $\Al(f)$.
For example, $3241$ is an allowed pattern of the {\em logistic map} $L(x)=4x(1-x)$ defined in $[0,1]$, obtained with initial value $x=0.8$.

\begin{figure}[htb]
\centering
\includegraphics[totalheight=52mm]{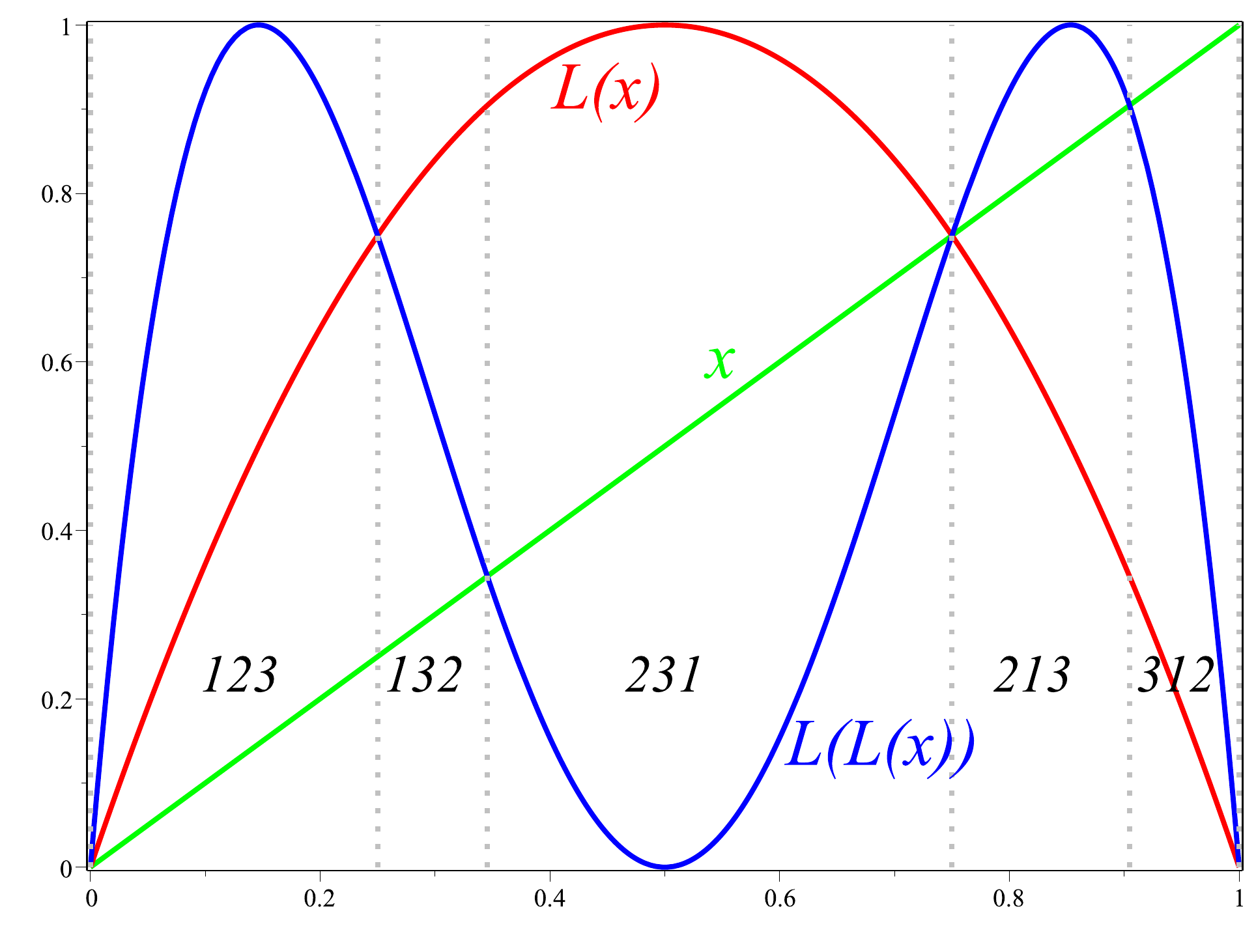}
\caption{The allowed patterns of length 3 of the map $L(x)=4x(1-x)$ on the unit interval.}
\label{fig:logistic}
\end{figure}

Bandt, Keller and Pompe~\cite{Bandt} proved that if $f$ is a piecewise monotone map on a one-dimensional interval, there are some permutations that are not realized by $f$. These permutations (i.e., those not in $\Al(f)$) are called {\em forbidden patterns}. For example, Figure~\ref{fig:logistic} shows that $321$ is a forbidden pattern of the logistic map $L(x)$. As shown in~\cite{AEK}, $L(x)$ and the so-called \emph{tent map} have the same allowed and forbidden patterns.

It is easy to see that, for every $f$, the set $\Al(f)$ is closed under consecutive pattern containment:  if $\pi\in\Al(f)$ and $\pi$ contains $\sigma$, then $\sigma\in\Al(f)$~\cite{Elishifts}. This property allows us to apply combinatorial tools to the study of dynamical systems. For example, it was shown in~\cite{Bandt} that for a piecewise monotone map $f$, the limit
\beq\label{eq:topo}
\lim_{n\rightarrow\infty} |\Al(f)\cap\S_n|^{1/n}
\eeq 
exists, and its logarithm equals the {\em topological entropy} of $f$, which is a measure of the complexity of the dynamical system. 
Additionally, the fact that $\Al(f)$ is closed under pattern containment allows us to characterize it by the set $\B(f)$ of minimal forbidden patterns, i.e., forbidden patterns where any proper pattern that they contain is allowed.
Thus, one can determine the topological entropy of a map by studying the asymptotic behavior of a set of (consecutive) pattern-avoiding permutations. Note that $\B(f)$ is an antichain in the consecutive pattern poset mentioned in Section~\ref{sec:other}. It was shown in~\cite{Eliu} that the logistic map defined above and, more generally, the family of maps $L_r(x)=rx(1-x)$, where $1<r\le4$, have infinitely many minimal forbidden patterns.

One of the goals of this combinatorial approach to dynamical systems is to obtain a better understanding of the sets of allowed and forbidden patterns of a map,
and in particular how the properties of these sets are related to the properties of the map.
Forbidden patterns are useful in certain tests to distinguish random from deterministic time series. These tests are based on the fact that in a sequence of values chosen independently at random from some continuous probability distribution, every permutation appears with positive probability. On the other hand, if the sequence has been generated by iterating a piecewise monotone map, then some patterns will never appear.
In fact, since $n!$ grows faster than $C^n$ for any constant $C$, the fact that the limit in Equation~\eqref{eq:topo} exists implies that most long patterns are forbidden.
Such permutation-based tests and their robustness against noisy data have been studied in~\cite{AZS,AZS2}, and they have been applied to the analysis of the stock market in~\cite{Zanin,Zanin2}.

Even though determining the allowed patterns of a map is a difficult problem in general, such patterns are relatively well understood for the so-called {\em shift maps}. The shift $\Sigma_N$ is defined on the set of infinite words over the alphabet $\{0,1,\dots,N{-}1\}$, ordered lexicographically, by $\Sigma_N(w_1w_2w_3\dots)=w_2w_3w_4\dots$. For example, $4217536\in\Al(\Sigma_3)$ because this permutation describes the relative order of the first few shifts of any word starting as $2102212210\dots$.
By considering words as expansions in base $N$ of real numbers in $[0,1]$, the shift $\Sigma_N$ is equivalent to the map that sends $x$ to the fractional part of $Nx$.
 It was shown in~\cite{AEK} that the shortest forbidden patterns of $\Sigma_N$ have length $N+2$, that there exactly six of this length, and that there are minimal forbidden patterns of each length $\ell\ge N+2$. For example, the patterns $1423$, $2134$, $2314$, $3241$, $3421$ and $4132$ are the shortest ones that cannot be realized by shifts of infinite binary words.
 In~\cite{Elishifts}, the allowed patterns of $\Sigma_N$ were characterized and enumerated, by means of a formula that determines the minimum number of symbols needed in the alphabet in order for a given permutation to be realized by a shift:

\begin{theorem}[\cite{Elishifts}]
For $\pi\in\S_n$, let $N(\pi)=\min\{k:\pi\in\Al(\Sigma_k)\}$. Then
$$N(\pi)=1+\des(w)+\epsilon,$$
where $w$ is the word $w_1\dots w_{\pi(n)-1}w_{\pi(n)+1}\dots w_n$ defined by $w_i=\pi(\pi^{-1}(i)+1)$ for all $i\neq\pi(n)$, and
$\epsilon$ equals $1$ if $\pi(n-1)\pi(n)\in\{21,(n{-}1)n\}$ and $0$ otherwise.
\end{theorem}

Replacing $N$ with an arbitrary real number $\beta>1$, one gets the {\em $\beta$-shift} $\Sigma_\beta$. Its combinatorial description is more elaborate~\cite{Par}, but its counterpart on the unit interval is simply the map $x \mapsto \fr{\beta x}$. In~\cite{Elibeta}, a method is given to compute, for any permutation~$\pi$, the smallest real number $\beta$ such that $\pi\in\Al(\Sigma_\beta)$, as a root of a polynomial whose coefficients depend on $\pi$. This value of $\beta$ is called the {\em shift-complexity} of $\pi$.
For example, the shift-complexity of $\pi=893146275$ is $\beta\approx3.343618091$, which is a root of $x^7-3x^6-x^5-2x^3+x^2+x-2$.

If instead of considering an arbitrary initial point in the domain of $f$ one restricts to periodic points, the permutations realized
by the relative order of the entries in the corresponding orbits (up until the first repetition) are called periodic patterns. 
For continuous maps, Sharkovskii's theorem~\cite{Sarko} describes the possible periods of these orbits, and there is significant literature that studies which periodic patterns are forced by others. From a more combinatorial perspective, Archer and Elizalde~\cite{ArcEli} characterize the periodic patterns of signed shifts, which are a large family of maps that includes both shifts and the tent map~\cite{Amigosigned}, 
in terms of the structure of the descent set of a certain cyclic permutation associated to the pattern. They also obtain enumeration formulas in some cases.

Another interesting direction of research considers, rather than iterates of a map $f$, random sequences given by a one-dimensional random walk that at time $i$ is at $X_1+X_2+\dots+X_i$, where the $X_j$ are i.i.d. continuous random variables. It is natural to ask what is the probability distribution on the permutations given by the relative order of the first $n$ values. In general, the probability of obtaining a given permutation depends on the distribution of the $X_j$, but there are pairs of permutations that are always equally likely to each other, such as $1432$ and $2143$. Elizalde and Martinez~\cite{EliM} define an equivalence relation on permutations under which equivalent patterns appear with equal frequency in any random walk, and 
characterize the resulting equivalence classes applying combinatorial methods.

We end by listing a few general questions in this area in the interface of permutation patterns and dynamical systems. They are mostly wide open, although some partial progress has been made in special cases~\cite{ABLLP,AEK,Elishifts,Elibeta,Eliu}:

\begin{enumerate}
\item Characterize the maps $f$ for which $\B(f)$ is finite.
\item Find a method to determine the length of the shortest forbidden pattern of a map. 
\item Characterize the sets of permutations that can be $\B(f)$ for some (piecewise monotone) $f$. 
\item Enumerate and/or characterize $\Al(f)$ or $\B(f)$ for particular maps~$f$.
\item Find the limit~\eqref{eq:topo} for particular piecewise monotone maps $f$, as a combinatorial way to determine their topological entropy.
\end{enumerate}

\begin{acknowledgement}
Partially supported by grant \#280575 from the Simons Foundation and by grant H98230-14-1-0125 from the NSA.
\end{acknowledgement}

\end{document}